\documentclass[12pt]{article}
\usepackage[cp1251]{inputenc}
\usepackage{amssymb}
\usepackage{amsthm}
\usepackage{inputenc}
\begin{document}
\newtheorem{defn}{Definition}[section]
\newtheorem{lm}{Lemma}[section]
\newtheorem{thm}{Theorem}[section]
\newtheorem{pr}{Proposition}[section]
\newtheorem{exam}{Example}[section]
\newtheorem{cor}{Corollary}[section]
\author{S. Albeverio $^{1},$ Sh. A. Ayupov $^{2},$ \ \ B. A.
Omirov  $^3$ \\ and R. M. Turdibaev$^4$}

\title{\bf CARTAN SUBALGEBRAS OF LEIBNIZ $n$-ALGEBRAS.}

\maketitle
\begin{abstract}
The present paper is devoted to the investigation of properties of Cartan subalgebras and regular elements in Leibniz $n$-algebras. The relationship between Cartan subalgebras and regular elements of given Leibniz $n$-algebra and Cartan subalgebras and regular elements of the corresponding factor $n$-Lie algebra is established.
\end{abstract}

\medskip
$^1$ Institut f\"{u}r Angewandte Mathematik, Universit\"{a}t Bonn, Wegelerstr. 6, D-53115 Bonn (Germany); SFB
611, BiBoS; IZKS, CERFIM (Locarno); Acc. Arch. (USI),  e-mail: \emph{albeverio@uni-bonn.de}

$^2$ Institute of Mathematics and Information Technologies, Uzbekistan Academy of Science, F. Hodjaev str. 29, 100125, Tashkent (Uzbekistan),
e-mail: \emph{sh\_ayupov@mail.ru}

 $^{3}$ Institute of Mathematics and Information Technologies, Uzbekistan
Academy of Science, F. Hodjaev str. 29, 100125, Tashkent (Uzbekistan), e-mail: \emph{omirovb@mail.ru}

 $^{4}$ Department of Mathematics, National University of Uzbekistan,
Vuzgorogok, 27, 100174, Tashkent (Uzbekistan), e-mail:
\emph{rustamtm@yahoo.com}



\section{Introduction}

The general notion of an algebra with a system $\Omega$ of polylinear 
operations was introduced by A.G. Kurosh in \cite{Kur} under the term of $\Omega$-algebra.
Starting with this notion V.T. Filippov \cite{Fil} defined an $n$-Lie algebra as an algebra with one $n$-ary polylinear operation ($n\geq 2$) which is anti-symmetric in all variables and satisfies a generalized Jacobi identity.In $n$-Lie algebras, similarly to the case of Lie algebras, operators of right multiplication are derivations with respect to the given $n$-ary multiplication and generate the Lie algebra of inner derivations. In this manner a natural generalization of Lie algebras was suggested for the case where the given multiplication is an $n$-ary operation. It worth mentioning that in the paper \cite{Nam} a natural example of an $n$-Lie algebra of infinitely differentiable functions in $n$ variables has been considered, in which the $n$-ary operation is given by the Jacobian. This $n$-Lie algebra is applied in the formalism of mechanics of Nambu, which generalizes the classical Hamiltonian formalism.

For further examples and methods of construction of $n$-Lie algebras we refer to \cite{Dal}-\cite{Fil}, \cite{Kas}.

The present work is devoted to a new algebraic notion -- so called Leibniz $n$-algebras which was introduced in \cite{Cas1} and has been further investigated in \cite{Cas2}-\cite{Cas3}, \cite{Poj}. These algebras are both "non antisymmetric" generalizations of $n$-Lie algebras and also of Leibniz algebras \cite{Lod}, which are determined by the following identity:

$$[x,[y,z]]=[[x,y],z] - [[x,z],y].$$

Investigations of Leibniz algebras and $n$-Lie algebras show that many properties of Cartan subalgebras and regular elements of Lie algebras may be extended to these more general algebras. Therefore a natural question occurs whether the corresponding classical results are valid for Cartan subalgebras and regular elements in Leibniz $n$-algebras. This problem is the main objective of this paper. In this direction, we establish the relationship between Cartan subalgebras and regular elements of given Leibniz $n$-algebra and Cartan subalgebras and regular elements of the corresponding factor algebra with respect to the ideal, generated by the elements of the form: $$[x_1, x_2, \dots, x_i,
\dots, x_j, \dots, x_n]+[x_1, x_2, \dots, x_j, \dots, x_i, \dots,
x_n].$$

All spaces and algebras in the present work are assumed to be finite dimensional.

\section{Preliminaries}

\begin{defn} {\em \cite{Cas1} A vector space $L$ over a filed $F$ with an $n$-ary multiplication $[-,-,...,-]: L^{\otimes n} \to L$ is called} a Leibniz $n$-algebra {\em if it satisfies the following identity}
$$ [[x_1,x_2, \dots ,x_n],y_2,\dots ,y_n]]=\sum_{i=1}^n [x_1,\dots ,
x_{i-1},[x_i,y_2,\dots,y_n],x_{i+1},\dots,x_n] \ \ \ (1)
$$
\end{defn}

It should be noted that if the product $[-,-,...,-]$ is antisymmetric in each pair of variables, i.e.
$$[x_1, x_2, \dots, x_i,
\dots, x_j, \dots, x_n]= - [x_1, x_2, \dots, x_j, \dots, x_i, \dots,
x_n].$$ then this Leibniz $n$-algebra becomes an $n$-Lie algebra.

\begin{exam}\label{exam0} {\em \cite{Cas1} Let $L$ be a Leibniz algebra with the product $[ - , - ].$ Then the vector space $L$ can be equipped with the Leibniz $n$-algebra structure with the following product:}
$$ [x_1, x_2, \dots, x_n]:= [x_1, [x_2, \dots, [x_{n-1},x_n]]].
$$
\end{exam}
Given an arbitrary Leibniz $n$-algebra $L$ consider the following sequences ($s$ is a fixed natural number, $1\leq s\leq n$):
$$ L^{<1>_s}=L,\qquad L^{<k+1>_s}=
\displaystyle[\underbrace{L,\ldots,L}_{(s-1)-times},L^{<k>_s},
\underbrace{L,\ldots,L}_{\scriptsize{(n-s)-\textrm{times}}}],
$$
$$L^1=L, \qquad L^{k+1}=\sum\limits_{i=1}^n\displaystyle[\underbrace{L,\ldots,L}_{\scriptsize{(i-1)-
\textrm{times}}},L^k,
\underbrace{L,\ldots,L}_{\scriptsize{(n-i)-\textrm{times}}}].$$

\begin{defn} {\em A Leibniz $n$-algebra $L$ is said to be} $s$-nilpotent {\em (respectively,} nilpotent{\em)  if there exists a natural number $k\in \mathbb{N}$ (respectively, $l\in \mathbb{N}$) such that $L^{<k>_s}=0$ (respectively, $L^l=0$).}
\end{defn}

It should be noted that for $n$-Lie algebras the above notions of $s$-nilpotency and nilpotency coincide. Let us also recall that for Leibniz algebras (i.e. Leibniz $2$-algebras) the notions of $1$-nilpotency and nilpotency were already known to coincide \cite{Ayup}.

The following example shows that the $s$-nilpotency property for Leibniz $n$-algebra ($n\geq 3$) essentially depends on $s.$

\begin{exam}\label{exam1} {\em Let $L_{p, q}$ be an $m$-dimensional algebra with the basis $\lbrace e_1,e_2,\\ \dots,e_m\rbrace$ ($m \geq n-1$) with the following product:
$$\begin{array}{lcr}
\ [e_{p-1},e_1, \dots ,e_{n-1}] & = & e_{p-1}\\
\ [e_{q-1},e_1, \dots ,e_{n-1}] & = & -e_{q-1},\\
\end{array}$$
where $2\leq p \ne q \leq n$  and all other products of basic elements assumed to be zero.}\end{exam} 

A straightforward calculation shows that this is a Leibniz $n$-algebra which is $s$-nilpotent for all $s \neq p, q$ but it is neither $p$-, nor $q$-nilpotent. Moreover $L_{p, q}$ is not nilpotent.

Everywhere below we shall consider only $1$-nilpotent $n$-algebras and thus for the sake of convenience "$1$-nilpotency"   will be called symply "nilpotentcy".

Set $A^{\times k}:=\underbrace{A\times A\times
\cdots \times A}_{k-{times}}.$

\begin{defn} {\em A linear map $d$ defined on a Leibniz $n$-algebra $L$ is called} a derivation {\em if 
$$d([x_1,x_2,\dots , x_n])=\sum_{i=1}^n [x_1,\dots d(x_i),\dots , x_n].$$
The space of all derivations of a given Leibniz $n$-algebra $L$ is denoted by $Der(L).$}
\end{defn}

The identity (1) and the properties of derivations easily imply the following equality:
$$[f,g]([x_1,...,x_n])=\displaystyle \sum_{i=1}^n [x_1,...,[f,g](x_i),..., x_n],\ \ \ \ (2)$$
for all $f,g \in Der(L),$ where $[f,g]=fg-gf$. Therefore the space $Der(L)$ forms a Lie algebra with respect to the commutator $[f,g]$.

Given an arbitrary element $x=(x_2,\dots,x_n) \in L^{\times(n-1)}$ consider the operator $R(x):
L \to L$ of right multiplication defined as $$R(x)(z)=[z,x_2,\dots,x_n].$$

Let $G$ be a subalgebra of the Leibniz $n$-algebra $L$, and put $R(G)=\{ R(g) \ | \ g \in G^{\times(n-1)}\}. $

The identity (1) implies that any right multiplication operator is a derivation and the identity (2) can be rewritten as 
$$[d,R(y_2,\dots, y_n)]=\sum_{i=2}^{n} R(y_2,\dots,d(y_i),\dots,
y_n). \ \ \ (3).$$

The space $R(L)$ is denoted further by $Inner(L)$ and its elements are called \textit{inner derivations}.

Further we have the following identity:
$$[R(x),R(y)]=-\sum_{i=2}^n R(x_2,...,R(y)(x_i),...,x_n) \ \ \ \ (4).$$

Thus $Inner(L)$ forms a Lie ideal in $Der(L).$

\begin{thm}[Engel's theorem] \label{thmEngel}  A Lie algebra $L$ is nilpotent if and only if each operator $R(x)$ is nilpotent for any $x \in L.$ \end{thm}

From the identity (4) and Theorem \ref{thmEngel} it follows that if $G$ is a nilpotent subalgebra of a Leibniz $n$-algebra $L$, then $R(G)$ is a nilpotent Lie algebra.

Indeed, from (4) we have 
$$[R(x),R(y)]=-\sum_{i_2+\dots+i_n=1}
R(R(y)^{i_2}(x_2),...,R(y)^{i_k}(x_k),...,R(y)^{i_n}(x_n)),$$
which implies easily that 
$$[[[R(a),\underbrace{R(b)],...,R(c)],R(d)]}_{\mbox{$m$-times}}=$$
$$ \begin{array}{rll}
(-1)^m & \cdot \displaystyle\sum_{i_2^{(1)}+\dots+i_n^{(1)}=1} &
\displaystyle\sum_{i_2^{(2)}+\dots+i_n^{(2)}=1} \cdots \\
  \cdots &\displaystyle\sum_{i_2^{(m)}+\dots+i_n^{(m)}=1}&
R((R(d)^{i_2^{(m)}}\circ R(c)^{i_2^{(m-1)}}\circ ...\circ
R(b)^{i_2^{(1)}} )(a_2),\dots\\
&&\dots ,(R(d)^{i_n^{(m)}}\circ R(c)^{i_n^{(m-1)}}\circ
...\circ R(b)^{i_n^{(1)}})(a_n)). \\
\end{array}$$

Therefore if $G^{\langle{p}\rangle_1}=0$ for an appropriate $p \in \mathbb{N},$ then for $m \ge (n-1)p$ each term on the right side becomes zero and thus $R(G)$ is nilpotent.

The following lemma gives a decomposition of a given vector space into the direct sum of two subspaces which are invariant with respect to a given linear transformation .

\begin{lm}[Fitting's Lemma] Let $V$ be a vector space and
$A:V\to V$ be a linear transformation. Then $V= V_{0A}\oplus V_{1A}$, where $A(V_{0A})\subseteq V_{0A},$
$A(V_{1A})\subseteq V_{1A}$ and $V_{0A}=\{ v\in V|\ A^i(v)=0$ for some $i\}$ and
$V_{1A}=\bigcap\limits_{i=1}^\infty A^i(V).$ Moreover, $A_{| V_{0A}}$ is a nilpotent transformation and
$A_{| V_{1A}}$ is an automorphism. \end{lm}
\begin{proof}
See \cite{Jac} (chapter II, \S 4).
\end{proof}

\begin{defn} {\em The spaces $V_{0A}$ and $V_{1A}$ are called
the} Fitting's null-com-\\ ponent {\em and the } Fitting's one-component {\em (respectively) of the space $V$ with respect to the transformation $A.$} \end{defn}

\begin{defn} {\em An element $h \in L^{\times (n-1)}$ is said to be} regular {\em for the algebra $L$ if the dimension of the Fitting's null-component of the space $L$ with respect to $R(h)$ is minimal.
In addition, its dimension is called a} rank {\em of the algebra $L.$}
\end{defn}

It is easy to see that the dimension of the Fitting's null-component of a linear transformation $A$ is equal to the order of the zero root of the characteristic polynomial of this transformation. Hence an element $h$ is regular if and only if the order of the zero characteristic root is minimal for $R(h).$

It should be noted that in the case of $n$-Lie algebras the operator $R(h)$ of right multiplication is degenerated. In particular, if the dimension of an $n$-Lie algebra $L$ is less than $n$ then we have $V_{0R(h)}=L.$
If $dim L \geq n$ then $\dim V_{0R(h)} \geq
n-1.$

Note also that for Leibniz algebras (i.e. $n=2$) the operator $R(h)$ is also degenerated \cite{Alb}.

Let us give an example of a Leibniz $n$-algebra ($n\geq 3$) which admits a non degenerated operator of right miltiplication.

\begin{exam}\label{exam2} {\em Consider an $m$-dimensional Leibniz $n$-algebra $L$ over a field $F$ with the following miltiplication:
$$[e_i,e_1,\dots , e_{n-1}]=\alpha_i e_i, \ \ \alpha\in F$$
where $\{e_1,\dots, e_m\}$ is the basis of the algebra, $\alpha_i \neq 0$ for all $1 \leq i \leq m,$ $\displaystyle
\sum_{i=1}^{n-1}\alpha_i=0$ and all other products are zero.

In this algebra the operator $R(e_1,\dots, e_{n-1})$ is nondegenerated.}
\end{exam}

We need the following lemma which can easily be proved.

\begin{lm} \label{lem2.2} Let $[-,-,\dots, -]: V^{\otimes n} \to V$  be a polylinear operation on a vector space $V$. The following conditions are equivalent:

$1) \ [x_1,\dots, x_i,  x_{i+1}, \dots, x_n]=-[x_1,\dots, x_{i+1}, x_i, \dots, x_n]$ \\for all $1 \leq i \leq
n-1$

$2) \ [x_1,\dots, x_i, \dots, x_j, \dots, x_n]=-[x_1,\dots, x_j, \dots, x_i, \dots, x_n]$ \\for all $1 \leq i
\neq j \leq n$

$3) \ [x_1,\dots, x_i, \dots, x_j, \dots, x_n]=0 $ if $x_i= x_j$ for some $1 \leq i \neq j \leq n$

$4) \ [x_1,\dots, x_i,  x_{i+1}, \dots, x_n]=0 $ if $x_i= x_{i+1}$ for some $1 \leq i \leq n-1.$

\end{lm}
 Consider the $n$-sided ideals
$$\begin{array}{rl}
 I = & ideal \langle
[x_1,\dots,x_i, \dots , x_j, \dots , x_n] \ | \ \exists i, j: x_i=x_j \rangle \\
J= & ideal \langle
[x_1,\dots,x_i, \dots , x_j, \dots , x_n] +[x_1,\dots,x_i, \dots , x_j, \dots,x_n ] \ | \  1\leq i \ne j \leq n \rangle \\
\end{array}$$
Lemma \ref{lem2.2} implies that these ideals coincide.
\begin{lm} \label{lem2.3} Let $L$ be a Leibniz $n$-algebra. If it admits a nondegenerated operator of right multiplication then $I=L$. \end{lm}
\begin{proof}
Let $x_2,...,x_n \in L$ be elements such that the operator $R(x_2,...,x_n)$ is non-degenerated.

Suppose first that $x_2,...,x_n$ are linearly dependent and $x_p=\displaystyle\sum_{i=2,\ i \ne p}^n \alpha_i x_i, \ \alpha_i\in \mathbb{C}.$ For any $a \in L $ there exists $b \in L$ such that 
$$a=[b,x_2,...,x_n]=\displaystyle \sum_{i=2,\, i\ne p}^n
\alpha_i[b,x_2,...,x_{p-1},x_i,x_{p+1},...,x_n] \in I$$ and therefore $L=I.$

Now suppose that $x_2,...,x_n$ are linearly independent. Since the operator $R(x_2,...,x_n)$ is nondegenerated, there exist $y_k \in L,$ $(2\leq k \leq n)$ such that $[y_k,x_2,..,x_n]=x_k.$

It is clear that the elements $y_2,...,y_n$ are also linearly independent.

Note that if $x_k \in I$ for some $k$, then $$L=[L,x_2,\dots,x_n] \subseteq I,$$ and therefore $L=I.$

Suppose that $x_k \not \in I$ for all $2\leq k \leq n.$ Consider the equalities 
$$x_k=[y_k,x_2,...,x_{k-1},x_k,x_{k+1},...,x_n]=\left[y_k,x_2,...,x_{k-1},[y_k,x_2,...,x_n],x_{k+1},...,x_n\right].$$
The identity (1) implies 
$$\left[ [y_k,x_2,\dots ,x_{k-1},y_k,x_{k+1},\dots ,x_n],x_2,...,x_n\right]=$$
$$=\left[[y_k,x_2,...,x_n],x_2,\dots ,x_{k-1},y_k,x_{k+1},\dots ,x_n]\right]+$$
$$+\left[y_k,[x_2,x_2,...,x_n],x_3\dots ,x_{k-1},y_k,x_{k+1},\dots ,x_n]\right]+\dots +$$
$$+\left[y_k,x_2,\dots ,x_{k-2},[x_{k-1},x_2,\dots ,x_n],y_k,\dots ,x_n\right]+$$
$$+\left[y_k,x_2,...,x_{k-1},[y_k,x_2,...,x_n],x_{k+1},...,x_n\right]+$$
$$+\left[y_k,x_2,\dots ,x_{k-1},y_k,[x_{k+1},x_2,\dots ,x_n],\dots ,x_n\right]+\dots +$$
$$+\left[y_k,x_2,\dots ,x_{k-1},y_k,x_{k+1},\dots ,x_{n-1},[x_n,x_2,\dots ,x_n]\right].$$

Since all summands on the right side except 
$$\left[ [y_k,x_2,\dots ,x_n],x_2,\dots, x_{k-1},y_k,x_{k+1},\dots ,x_n\right],$$
$$\left[y_k,x_2,...,x_{k-1},[y_k,x_2,...,x_n],x_{k+1},...,x_n\right]$$
and the whole left side in the above equality belong to the ideal $I$, it follows that $$\left[y_k,x_2,...,x_{k-1},[y_k,x_2,...,x_n],x_{k+1},...,x_n\right]+$$
$$+\left[[y_k,x_2,...,x_n],x_2,\dots ,x_{k-1},y_k,x_{k+1},\dots ,x_n\right]\in
I.$$
Therefore $$x_k+[x_k,x_2,...,x_{k-1},y_k,x_{k+1},...,x_n]\in I. \ \ \ \ (5)$$

If $a \in \langle {x_2,...,x_n}\rangle \cap I$ then $a=\displaystyle \sum_{i=2}^n\alpha_ix_i.$

We have that $$[a,x_2,\dots ,x_{k-1},y_k,x_{k+1},\dots ,x_n] \in I.$$ 
On the other hand $$[a,x_2,\dots , x_{k-1},y_k,x_{k+1},\dots ,x_n]=\sum_{i=2}^n \alpha_i[x_i,\dots ,x_{k-1},y_k,x_{k+1},\dots ,x_n].$$
Therefore
$$\alpha_k[x_k,\dots ,x_{k-1},y_k,x_{k+1},\dots ,x_n] \in I.$$
Thus (5) implies that $\alpha_k x_k \in I$ and hence 
$\alpha_k=0 \ (2 \leq k \leq n),$ i.e. $a=0.$ This means that $$\langle {x_2,...,x_n}\rangle \cap I=0.$$

Note that $y_k \not \in I,$ because if $y_k \in I$ then (5) implies that $x_k \in I$ which contradicts our assumption that $x_k \notin I$ for all $k.$

If $a \in \langle  y_2,\dots , y_n \rangle \cap \langle
x_2,\dots , x_n \rangle,$ then $a=\displaystyle
\sum_{i=2}^n\beta_iy_i$ and $a=\displaystyle\sum_{i=2}^n\alpha_ix_i$
for some $\alpha_i, \beta_i \in$ $\mathbb{C}.$

By applying the operator $R(x_2,...,x_n)$ to the element $a$ we obtain $$\displaystyle \sum_{i=2}^n\beta_i[y_i,x_2,\dots
,x_n]=\sum_{i=2}^n\alpha_i[x_i,x_2,...,x_n] \in I,$$ i.e. $\displaystyle\sum_{i=2}^n\beta_ix_i \in I.$ But 
$\langle {x_2,...,x_n}\rangle \cap I=\lbrace0\rbrace$ and therefore $\beta_i=0$ for all $ 2\leq i \leq n$ and thus 
$$\langle  y_2,\dots , y_n \rangle \cap \langle x_2,\dots , x_n
\rangle =0.$$

If $a \in \langle {y_2,...,y_n}\rangle \cap I$ then 
$a=\displaystyle\sum_{i=2}^n\alpha_iy_i$ and $$\sum_{i=2}^n
\alpha_ix_i =\sum_{i=2}^n \alpha_i[y_i,x_2,\dots ,x_n]=
[a,x_2,\dots ,x_n] \in I.$$ Therefore $\alpha_i=0,$ i.e.
$a=0$ and $$\langle {y_2,...,y_n}\rangle \cap I=0.$$

Since $R(x_2,...,x_n)$ is nondegenerated, there exist $z_k \in L, (2\leq k\leq n)$ such that $[z_k,x_2,..,x_n]=y_k.$

In the same way one can prove that $z_2,\dots , z_n$ are linearly independent and 
$$ \langle  z_2,\dots , z_n \rangle \cap \langle  x_2,\dots , x_n \rangle=0,$$
$$ \langle  z_2,\dots , z_n \rangle \cap \langle  y_2,\dots , y_n \rangle=0,$$
$$ \langle  z_2,\dots , z_n \rangle \cap I=0.$$

Repeating the above process we obtain that $$I \oplus \langle x_2,\dots ,
x_n \rangle \oplus \langle y_2,\dots , y_n \rangle \oplus \langle
z_2,\dots , z_n \rangle \oplus \cdots \subseteq L,$$ which contradicts the finiteness of  $dim L.$ The proof is complete.
\end{proof}
The following example shows that the converse assertion to the Lemma \ref{lem2.3} is not true in general.

\begin{exam} {\em Consider a complex $m$-dimensional ($m\geq 4$) non Lie Leibniz algebra $L$ with the basis $\lbrace e, f, h, i_0, i_1, \dots, i_{m-4}\rbrace$ and the following table of multiplication:
$$[i_k, h]=(m-4-2k)i_k, \ \ 0 \leq k \leq m-4;$$
$$[i_k, f]=i_{k+1}, \ \ 0 \leq k \leq m-5;$$
$$[i_k, e]=k(k+3-n)i_{k-1}, \ \ 1 \leq k \leq m-4;$$
$$[e, h]=2e, \ \ [h, e]=-2e, \ \ [f, h]= -2f,$$
$$[h, f]=2f, \ \ [e, f]=h, \ \ [f, e]= -h,$$
with other products are zero.

Following the construction of Leibniz $n$-algebras from Example \ref{exam0} we obtain a Leibniz $n$-algebra. For $n> 4$ one has

$$I \ni [h,h,\dots,h,e]=[h,[h,\dots,[h,e]\dots]]=(-2)^{n-1}e,$$
$$I \ni [h,h,\dots,h,f]=[h,[h,\dots,[h,f]\dots]]=2^{n-1}f,$$
$$I \ni [f,h,\dots,h,e]=[f,[h,\dots,[h,e]\dots]]=-(-2)^{n-2}h,$$
$$I \ni [i_k ,h,\dots,h,f]=[i_k,[h,\dots,[h,f]\dots]]=2^{n-2}i_{k+1}$$ for $0\leq k \leq m-5$ and $$I \ni [i_1 ,h,\dots,h,e]=[i_1,[h,\dots,[h,e]\dots]]=(-2)^{n-2}(4-n)i_0.$$

Therefore for this Leibniz $n$-algebra ($n>4$) we have $I=L.$} \end{exam}

Moreover, let us show that in this Leibniz $n$-algebra all operators of right multiplication are degenerated.
Indeed, suppose that for some $a=(a_2,\dots,a_n) \in L^{\times (n-1)}$ the operator $R(a)$ is nondegenerated. Then for every $x\in L$ we have $$0\neq [x,a_2,\dots,a_n]=[x,[a_2,\dots,[a_{n-1},a_n]\dots]],$$
and for the element $b=[a_2,\dots,[a_{n-1},a_n]\dots]$ we obtain that $R(b)$ is a non degenerated operator in the Leibniz algebra $L$, which contradicts the Lemma 2.6 from \cite{Alb}.

We also have the following generalization of Fitting's Lemma for Lie algebras of nilpotent transformations of a
vector space.

\begin{thm} \label{thm2.2} Let $G$ be a nilpotent Lie algebra of linear
transformations of a vector space $V$ and $V_0=\bigcap\limits_{A \in G} V_{0A},$ $V_1=\bigcap\limits_{
i=1}^\infty G^i(V).$ Then the subspaces $V_0$ and $V_1$ are invariant with respect to $G$ (i.e. they are
invariant with respect to every transformation $B$ of $G$) and $V=V_0\oplus V_1.$ Moreover,
$V_1=\sum\limits_{A\in G} V_{1A}.$ \end{thm}
\begin{proof}
See \cite{Jac} (chapter II, \S 4).
\end{proof}

\emph{Remark 1.} From \cite{Jac} (chapter III, p. 117) in the case of a vector space V over an
infinite field and under the conditions of Theorem
\ref{thm2.2}, we have the existence of an element $B\in G$ such
that $V_0= V_{0B}$ and $V_1=V_{1B}$.

\section{The main results.}
\vspace{1cm}

Let $\Im$ be a nilpotent subalgebra of an $n$-Leibniz algebra $L$
and $L=L_0\oplus L_1$ be the Fitting's decomposition of $L$
with respect to the nilpotent Lie algebra $R(\Im)=\{ R(x)|\ x\in
\Im^{\times(n-1)} \}$ of transformations of the underlying vector space $V$ as in
Theorem \ref{thm2.2}.

\begin{defn}  {\em Given a subset $X$ in a Leibniz $n$-algebra $L,$ the} $s$-normalizer {\em of $X$ is the set $$ N_s(X)=\{a\in L \, |\, [x_1,\dots,x_{s-1},a,x_{s+1},\dots,x_n]\in X \textrm{for all} \ x_i \in X\}.$$} \end{defn}

Note that if in the Example \ref{exam2} we consider the set $X$ 
generated by the vectors $\langle e_1,e_2,
\dots,e_{n-1}\rangle,$ then $ N_1(X)=X$ and $N_s(X)=L$ for all $1 < s
\leq n.$

Further we shall consider only $1$-normalizers, and therefore we shall call them simply normalizers and denote the set of normalizer by $N(X).$

\begin{defn} {\em A subalgebra $\Im$ of a Leibniz $n$-algebra $L$ is said to be a} Cartan subalgebra {\em if :\\
\indent a) $\Im$ is nilpotent;\\
\indent b) $\Im=N(\Im).$} \end{defn}

The following example shows the existence of such subalgebras.
\begin{exam} {\em Consider the algebra $L=\langle e_1,e_2,\dots,e_m \rangle$ with the following multiplication:
$$[e_k,e_1,e_1,\dots,e_1]=e_k \ \ \ (2 \leq k \leq m).$$ It easy to see that $L$ is neither a nilpotent Leibniz $n$-algebra and nor an $n$-Lie algebra.

Consider the subspace $H=\langle e_1 \rangle.$ It is clear that $H$ is a nilpotent subalgebra. Put $$a=\alpha e_1 +\sum_{k=2}^m \beta_ke_k\in N(H),$$ then $H\ni [a,e_1,\dots,e_1]=\displaystyle\sum_{k=2}^m \beta_ke_k$ and therefore  $\beta_k=0$ for all $2\leq k \leq
m.$ Thus $H=N(H)$ and $H$ is a Cartan subalgebra in $L.$}
\end{exam}

\begin{lm} Let $L$ be the Leibniz $n$-algebra ($n\geq 3$) constructed from a Leibniz algebra as in Example \ref{exam0} and let $\Im$ be a Cartan subalgebra of the Leibniz algebra $L.$ Then in the Leibniz $n$-algebra we have:

a) $\Im$ is a nilpotent subalgebra;

b) $N(\Im)=L.$
\end{lm}
\begin{proof} The nilpotency of the subalgebra $\Im$ in the Leibniz $n$-algebra follows from its nilpotency in the Leibniz algebra $L$. From \cite{Omi} it is known that under the natural homomorphism of a Leibniz algebra onto the corresponding factor algebra which is a Lie algebra, the image of the Cartan subalgebra $\Im$ is a Cartan subalgebra of the Lie algebra. Further using abelianness of Cartan subalgebras in Lie algebras \cite{Jac} we obtain that $[\Im,\Im]$ is contained in the ideal generated by the squares of elements from the Leibniz algebra $L,$ and this ideal is contained in the right annihilator. Therefore for $n\geq 3$ we have $N(\Im)=\{x\in L | [x, \Im, \dots, \Im]\subseteq
\Im \}=L.$
\end{proof}

For Cartan subalgebras of $n$-Leibniz algebras similar to the case
of $n$-Lie algebras and Leibniz algebras, there is a
characterization in terms of the Fitting's null-component, namely,
the following proposition is true.
\begin{pr} \label{prop3.1} Let $\Im$ be a nilpotent subalgebra of a Leibniz $n$-algebra $L.$ Then $\Im$ is a Cartan subalgebra if and only if it coincides with $L_0$ in the Fitting decomposition of the algebra $L$ with respect to $R(\Im).$
\end{pr}

\begin{proof}
Let $x \in N(\Im),$ then $[x,h_2,...,h_n] \in \Im$ for all 
$h_i \in \Im \ (2 \leq i \leq n).$  Since $\Im$ is nilpotent there exists $k \in \mathbb{N}$ such that  $R^k(h_2,...,h_n)(x)=0,$ i.e. $x \in L_0.$

Therefore we have $N(\Im) \subseteq L_0.$ Since $\Im \subseteq N(\Im)$ we obtain that $\Im \subseteq
L_0.$

Suppose that $\Im \subsetneqq L_0.$

Taking $R(\Im)$ instead of $G$ in Theorem \ref{thm2.2} we obtain that $L_0$ is invariant with respect to $R(\Im)$ and $R(h_2,...,h_n)|_{L_0}$ is a nilpotent operator for all $h_i \in \Im \ (2\leq i \leq n).$

Therefore we have $$R(\Im): L_0 \to L_0,$$ $$R(\Im): \Im
\to \Im,$$ where $R(\Im)$ is a Lie algebra.

Thus we obtain the induced Lie algebra $\overline {R(\Im)}:
L_0/\Im \to L_0/\Im ,$ where $L_0/\Im$ is a non-zero linear factor space. If we consider $R(\Im): L_0 \to L_0,$ then as it was mentioned above the operator $R(h_2,...,h_n)$ is nilpotent for all $h_i\in \Im \ (2\leq i \leq n).$ Then by Engel theorem \cite{Jac} it follows that there exists a non-zero element $\overline
x=x+\Im$  $(x\notin \Im)$ such that $\overline{R(\Im)}(x+\Im)=\overline{0}.$ This means that $[x,h_2,...,h_n] \in \Im$ for every $h_i \in \Im \ (2\leq i \leq n).$
Therefore there exists an element $x \in N(\Im)$ such that $x \not\in \Im$ -- the contradiction shows that $\Im =L_0.$ The proof is complete. \end{proof}

\begin{cor} \label{cor11}Let $\Im$ be a Cartan subalgebra of the Leibniz $n$-algebra $L$. Then $\Im$ is a maximal nilpotent subalgebra of $L$.
\end{cor}
\begin{proof} Let $B$ be a nilpotent subalgebra of the $L$ such that $\Im \subseteq B$ then by Proposition \ref{prop3.1} we have $\Im \subseteq B \subseteq L_0(\Im)=\Im.$ 
\end{proof}
The following theorem establishes properties of the Fitting's null-component of the regular element of an $n$-Leibniz algebra.

\begin{thm} Let $L$ be a Leibniz $n$-algebra over an infinite field and let $x$ be a regular element for $L.$ Then the Fitting null-component $\Im=L_0$ with respect to the operator $R(x)$ is a nilpotent subalgebra in $L.$
\end{thm}

\begin{proof}
Let us prove that both Fitting components with respect to $R(x)$ are invariant under $R(\Im).$ Indeed, let $a=(a_2,...,a_n) \in \Im^{\times (n-1)}.$ Then from (4) it easily follows that  $$[[[R(a),\underbrace{R(x)],R(x)],...,R(x)}_{m-times}]= (-1)^m\displaystyle\sum_{i_2+\dots+i_n=m}
R\left(R(x)^{i_2}(a_2),...,R(x)^{i_n}(a_n)\right).$$

For sufficiently large $m$ we obtain that $$[[[R(a),\underbrace{R(x)],R(x)],...,R(x)}_{ m-{times}}]=0.$$

From \cite{Jac} (Lemma 1, chapter II, \S 4) we have that the Fitting components $L_0$ and $L_1$ with respect to $R(x)$ are invariant under $R(a).$

Let us prove that the operator $R(h_2,...,h_n)|_{L_0}$ is nilpotent for $h_i \in \Im \
(2\leq i \leq n)$. Assume the opposite, i.e. there exists $h=(h_2,...,h_n)\in \Im^{\times (n-1)}$ such that $R(h_2,...,h_n)|_{L_0}$ is not nilpotent.

Consider $u^t=(u^t _2,...,u^t _n),$ where $u^t _i=tx_i+(1-t)h_i$ and $t$ belongs to the underlying field. Then the elements of the matrices $R(u^t _2,...,u^t _n)|_{L_0}$ and $R(u^t _2,...,u^t _n)|_{L_1}$
are polynomials in $t.$ Since  for $t=1$ we have $R(x)|_{L_1}=R(u^1)|_{L_1}$ and $R(u^1)|_{L_1}$ is non degenerated, there exists $t_0$ such that $R(u^{t_0})|_{L_1}$ is non degenerated and $R(u^{t_0})|_{L_0}$ is not nilpotent.

In this case the dimension of the Fitting null-component of the space $L$
with respect to the operator $R(u^{t_0})$ is less than the dimension of the Fitting null-component with respect to the operator $R(a),$ which contradicts the regularity of the element $x.$

Therefore, $R(h)|_{L_0}$ is nilpotent for every $h\in
\Im^{\times (n-1)}$ and by Theorem \ref{thmEngel} $\Im$ is nilpotent. The proof is complete.\end{proof}

Let us recall that the Fitting null-component with respect to the right multiplication operator by a regular element in $n$-Lie algebras \cite{Kas} and Leibniz algebras \cite{Alb} is a Cartan subalgebra. But in the case of Leibniz $n$-algebras the Example \ref{exam2} shows that the Fitting null-component with respect to the operator of right multiplication by the regular element $e=(e_1, e_2, \dots, e_{n-1})$ is not a Cartan subalgebra, because $V_{0R(x)}=\{0\}$ and $N(\lbrace0\rbrace)=L$.

\begin{pr}\label{prop3.2} Let $L$ be a Leibniz $n$-algebra over a field $F$ and let $\Omega$ be an arbitrary extension of the field $F$. Put $L_{\Omega}=L_F\otimes \Omega.$ Then $\Im$ is a Cartan subalgebra in $L$ if and only if $\Im_{\Omega}=\Im_F \otimes \Omega$ is a Cartan subalgebra in $L_{\Omega}$.\end{pr}

\begin{proof} Let $\Im$ be a Cartan subalgebra in $L.$ Then $\Im_\Omega$ is a subalgebra in $ L_\Omega.$

Since $\Im^{<k>_1}=0,$  from the evident equality $\Im^{<k>_1}_\Omega=\Im^{<k>_1}_F \otimes \Omega $ it follows that $\Im_\Omega$ is a nilpotent subalgebra.

Consider $a_\Omega \in N(\Im_\Omega).$ Then $[a_\Omega, h_{2\Omega },\dots,h_{n\Omega }]=[a, h_2,\dots,h_n]\otimes \alpha\gamma_2\dots\gamma_n,$ where $a_\Omega=a \otimes \alpha,\,
h_{i\Omega}=h_i\otimes \gamma_i \ (2 \leq i\leq n).$
Therefore, $ N(\Im_\Omega ) \subseteq N(\Im) \otimes \Omega.$
But $\Im$ is a Cartan subalgebra and thus $N(\Im)=\Im$ and $N(\Im_\Omega ) \subseteq \Im_F\otimes \Omega=\Im_\Omega.$

Therefore, $\Im_\Omega $ is a Cartan subalgebra in $L_{\Omega}.$

Conversely, suppose that $\Im_\Omega$ is a Cartan subalgebra in $L_{\Omega}.$ Then from $\Im^{<k>_1}_F \otimes \Omega=\Im^{<k>_1}_\Omega=0$ it follows that $\Im^{<k>_1}_F=0.$

Consider $a \in N(\Im_F).$ We have $[a,h_2,\dots,h_n]\in \Im_F$
for all $h_i \in \Im \ (2\leq i \leq n).$ Since $[a_\Omega, h_{2\Omega },\dots,h_{n\Omega }]=[a,
h_2,\dots,h_n]\otimes \alpha\gamma_2\dots\gamma_n \in \Im_F
\otimes \Omega,$ where $a_\Omega=a \otimes \alpha,\,
h_{i\Omega}=h_i\otimes \gamma_i \ (2 \leq i\leq n),$ one has $a\otimes
\alpha = a_\Omega \in N(\Im_\Omega)=\Im_\Omega=\Im_F \otimes
\Omega.$ Therefore, $a \in \Im_F$ and $\Im_F$ is a Cartan subalgebra in $L_F.$ The proof is complete. \end{proof}

\begin{thm} \label{thm3.2}Let $\varphi: L \to L'$ be an epimorphism of Leibniz $n$-algebras and suppose that $\Im$ is Cartan subalgebra in $L$ and $\varphi(\Im)=\Im'.$ Then $\Im'$ is a Cartan subalgebra in $L'.$ \end{thm}
\begin{proof} In view of Proposition \ref{prop3.2} we may assume that the field $F$ is algebraically closed.

Consider the decomposition into the sum of characteristic subspaces:
$$L=L_\alpha \oplus L_\beta \oplus \cdots \oplus L_\gamma$$
with respect to the nilpotent Lie algebra $R(\Im)$ of linear transformations of the vector space $L,$ where $L_\alpha=\{x \in L \, | \, \left(R(h)-
\alpha(h)id\right)^k(x)=0$ for some $k$ and for any $h \in \Im^{\times (n-1)}\}.$ Then $$\varphi(L)=\varphi(L_\alpha)+\varphi(L_\beta)+\cdots+\varphi(L_\gamma).$$

By using the properties of homomorphisms we obtain by induction that $$\varphi \circ R(x_2,...,x_n)^k= R(\varphi(x_2),...,\varphi(x_n))^k \circ \varphi$$ for every $k \in
\mathbb{N}.$

Further we have $$\begin{array}{rcl}
& \varphi \circ (R(h)-\alpha id)^k =   & \varphi \circ \displaystyle\sum_{i=0}^k C_k^i\alpha^{n-k}R(h)^k=\\
&\displaystyle\sum_{i=0}^k C_k^i\alpha^{n-k} \varphi \circ R(h)^k = & \displaystyle\sum_{i=0}^k C_k^i\alpha^{n-k}R(\varphi(h))^k\circ \varphi=\\
&\displaystyle(R(\varphi(h))-\alpha id)^k\circ \varphi, & \\
\end{array}$$
where $C_k^i$ are binomial coefficients.

Therefore from $$ (R(h)-\alpha(h)id)^k(x)=0$$ we obtain 
$$(R(h')-\alpha(h')id)^k\varphi(x)=0,$$ where $h'=\varphi(h).$

Thus, if $x\in L_\alpha,$ then $x'\in L_\alpha '$ (where $\varphi(L_\alpha)=L_\alpha'$). Since $\varphi$ is epimorphic, we have the following decomposition of the space $L'$ with respect to  $R(\Im')$:
$$L'=L'_\alpha \oplus L'_\beta \oplus \cdots \oplus L'_\gamma,$$ where $\varphi(L_\sigma)=L'_\sigma$ and $\sigma\in \{\alpha, \beta, \dots,
\gamma\}.$

If $\alpha \ne 0,$ the action of $\Im'$ on $L'_\alpha$ is non degenerated and therefore $L_0'=\varphi(L_0)=\varphi(\Im)=\Im'.$
Now Proposition \ref{prop3.1} implies that $\Im'$ is a Cartan subalgebra of $L'$. The proof is complete. \end{proof}

\begin{pr} \label{prop3.3} If the sum $x_1+x_2+\ldots+x_m=x$ of characteristic vectors corresponding to different characteristic values $\rho_1,\rho_2,\ldots,\rho_m$ of a transformation $Q$ belongs to an invariant subspace $M,$ then each summand is contained in $M.$\end{pr}
\begin{proof} See \cite{Mal}, Chapter III, p. 147.\end{proof}

\begin{pr} \label{prop3.4} Let $M$ be an invariant subspace of a vector space $V$ with respect to a linear transformation $Q:V \to V.$ If for the element $x=x_0+x_\alpha+x_\beta+\cdots+x_\gamma$ decomposed into the sum of characteristic vectors from corresponding characteristic spaces $V_\zeta$ ($\zeta \in
\{0,\alpha,\beta,...,\gamma\}$) we have $Q(x)\in M.$ Then $x-x_0 \in M.$ \end{pr}
\begin{proof} If $x=x_0$ then everything is clear.

Suppose that $x \neq x_0,$ i.e. there exists a non zero vector among $x_\alpha,x_\beta,  \dots, \\ x_\gamma.$  We may assume that $x_\alpha \neq 0.$

From $Q(x)\in M$ it follows that $$Q(x_0)+Q(x_\alpha)+Q(x_\beta)+\ldots+Q(x_\gamma)\in M.$$

Since $M$ is invariant and $Q(x_\sigma) \in L_\sigma$ for every $\sigma \in \{0,\alpha,\beta,...,\gamma\},$  Proposition \ref{prop3.3} implies that $Q(x_\sigma) \in M$ for every 
$\sigma \in \{0, \alpha,\beta,...,\gamma\}.$

Note that $Q(x_\alpha) \neq 0$ because $\alpha=0.$

We have $Q(M\cap L_\alpha) \subseteq Q(M)\cap Q(L_\alpha) \subseteq
M \cap L_\alpha.$ It is clear that for non zero $y \in M\cap
L_\alpha$ we have $Q(y)\neq 0$ (because  otherwise $y \in L_0,$ while $L_0\cap
L_\alpha=0$).

Therefore $Q(M \cap L_\alpha) = M \cap L_\alpha.$ Since $Q(x_\alpha) \in M \cap L_\alpha$ there exists $y\in M \cap L_\alpha$ such that $$Q(x_\alpha)=Q(y) \Rightarrow
Q(x_\alpha-y)=0  \Rightarrow x_\alpha-y \in L_0.$$

But since $x_\alpha, y \in L_\alpha,$ we have $x_\alpha-y \in L_\alpha \cap L_0.$ Therefore, $x_\alpha=y \in M.$

Since $\alpha$  is arbitrary we obtain that $x_\alpha,x_\beta, \dots ,x_\gamma \in M$ and therefore $x-x_0=x_\alpha+x_\beta+\ldots x_\gamma \in M.$ The proof is complete. \end{proof}

For a Leibniz $n$-algebra $L$ consider the natural homomorphism $\varphi$ onto the factor algebra $\overline{L}=L/I.$ It is clear that $\overline{L}$ is an $n$-Lie algebra.

\begin{cor} \label{cor3.1} Let $b \in L^{\times(n-1)}$.  Consider the decomposition of the element $x=x_0+x_\alpha+x_\beta+\cdots +x_\gamma$ with respect to $R(b),$ where  $x_\sigma\in L_\sigma,$ $\sigma \in \{0,\alpha,\beta,...,\gamma\}$. If there exists $k \in \mathbb{N}$ such that $R(b)^k(x) \in I,$ then $\overline x=\overline x_0.$
\end{cor} 
\begin{proof} Let $R(b)^k(x) \in I$ and $R(b)^{k-1}(x) \notin I.$

Setting $Q:=R(b)^k,$ we obtain $Q(x)\in I.$ On the other hand $Q(I)\subseteq I$ since $I$ is an ideal in $L.$ Proposition \ref{prop3.4} implies that $x-x_0 \in I,$ i.e. $\overline x=
\overline x_0.$ The proof is complete. \end{proof}

\emph{Remark 2.} For the Cartan subalgebra $\Im$ of the
Leibniz $n$-algebra $L,$ we consider the Lie algebra $R(\Im)$ of
linear transformations $L$ (which evidently is nilpotent) and the
decomposition of $L$ with respect to $R(\Im).$ Remark 1
implies the existence of an element $R(b)\in R(\Im)$ such that the
Fitting's null-component with respect to the nilpotent Lie algebra
of linear transformations $R(\Im)$ coincides with the Fitting's
null component with respect to the transformation $R(b),$ i.e.
$L_0= L_0(b).$ Using Proposition \ref{prop3.1} we obtain $\Im=
L_0(b).$

\begin{lm} \label{lm3.2} Let $b \in \Im^{\times(n-1)}$ and $\Im=L_0(b).$ Then $\overline \Im=L_{\overline0}(\overline b).$ \end{lm}
\begin{proof} Let $\overline \Im$ be the image of the Cartan subalgebra $\Im$ under the homomorphism $\varphi :L \to L/I.$ From the theory of $n$-Lie algebras \cite{Kas} we know that there exists a regular element $\overline a=(\overline a_2, \overline a_3, \dots, \overline a_n)  \in
\overline \Im^{\times(n-1)}$ such that $\overline
\Im=L_{\overline0}(\overline a).$

Without loss of generality we may assume that $a=(a_2, a_3, \dots,
a_n) \in \Im^{\times(n-1)}.$ It is clear that $L_0(b)\subseteq
L_0(a).$ Since $\overline a$ is a regular element we have that $\overline
\Im \subseteq L_{\overline0}(\overline b).$

If there exists $i$ such that $a_i \in I,$  then $\overline
L=L_{\overline0}(\overline a)\subseteq L_{\overline0}(\overline
b)$ and therefore $L_{\overline0}(\overline
a)=L_{\overline0}(\overline b)$ and $\overline{\Im}=L_{\overline0}(\overline b).$

Suppose that for any $i$ we have $a_i \not \in I$ and $\overline{\Im}\subsetneqq L_{\overline0}(\overline b)$. Then there exists $x$ such that $\overline x=x+I \in
L_{\overline0}(\overline b) \setminus L_{\overline0}(\overline
a).$ Therefore for the element $x$ we have $R(b)^k(x) \in I$ for some $k$ and $R(a)^s(x) \not \in I$ for any $s\in \mathbb{N}.$

Note that $R(b)^t(x) \neq 0$  for any $t \in \mathbb{N},$ because in the other case $x \in L_0(b)\subseteq L_0(a)$ which contradicts the condition $\overline x \not \in L_{\overline0}(\overline a).$ Therefore $x \notin \Im.$

Thus for the element $x$ we have $R(b)^k(x) \in I$ and $x\ne x_0.$ Corollary \ref{cor3.1} implies that $\overline x=\overline x_0\in \overline \Im =L_{\overline0}(\overline a)$,
which contradicts the choice of $x.$ Therefore $\overline \Im=
L_{\overline0}(\overline b).$ The proof is complete. \end{proof}

\begin{thm} The image of a regular element for a Leibniz $n$-algebra $L$ under the natural homomorphism $\varphi : L \rightarrow L/I$ is a regular element for the $n$-Lie algebra $L/I.$ \end{thm}
\begin{proof} Suppose that $a=(a_2, a_3, \dots, a_n) \in L^{\times (n-1)}$ is a regular element for $L$ and $\overline a=(a_2+I, a_3+I, \dots, a_n+I)$ is not regular in $L/I.$ Let $\overline b=(b_2+I, b_3+I, \dots, b_n+I)$ be an arbitrary regular element for $L/I,$ then $a_i-b_i \notin I$ for some $i$.

Since $I$ is an ideal in $L$, for every $x\in L^{\times (n-1)}$ we have $R(x)(I) \subseteq I$ and the matrix of the operator $R(x)$ in the basis $\{ e_1,e_2,\ldots, e_m,i_1,i_2,\ldots, i_l\}$ of the algebra $L$ (where $\{ i_1,i_2,\ldots, i_l\}$ the basis of $I$) has the following form:
$$ R(x)=
\left(\begin{array}{cc}
X,& 0 \\
Z_x, & I_x \\
\end{array} \right), $$
where $X$ is the matrix of the operator $R(x)|_{\{e_1,\dots,e_m\}}$ and $I_x$ is the matrix of the operator $R(x)|_I.$

Let $$ R(a)=\left(\begin{array}{cc} A,& 0 \\ Z_a, & I_a \\
\end{array}
\right),\quad R(b)=\left(\begin{array}{cc} B,& 0 \\ Z_b, & I_b \\
\end{array} \right)
$$ be the matrices of the transformations $R(a)$ and $R(b)$ respectively.

Denote by $k$ (respectively by $k'$) the order of the $0$ characteristic value  of the matrix $A$ (respectively $B$) and by $s$ (respectively by $s'$) the order of the $0$ characteristic value of the matrix $I_a$ (respectively $I_b$). Then we have $k'<k, \ s<s'.$

Put $U=\bigg\{ y\in L^{\times (n-1)}\left|\ R(y)=\left(\begin{array}{cc} Y,& 0 \\
Z_y, & I_y \\ \end{array} \right)\right.$ and $Y$ has the $0$ characteristic value of the order less than $k \bigg\}$ and $ V=\bigg\{y\in L^{\times (n-1)} \left| \ R(y)=\left(\begin{array}{cc} Y,& 0 \\ Z_y, & I_y \\ \end{array} \right)\right.$ and $I_y$ has the $0$ characteristic value of the order less than  $s+1 \bigg\}.$

Since $b\in U$ and $a\in V,$ the above sets are non empty. Similar to considerations in \cite{Alb} one can prove that the sets $U$ and  $V$ are open in the Zariski topology and therefore they have non-empty intersection. Let $y\in U\cap V,$ i.e. $y\in L^{\times (n-1)}$ is such an element that $Y$ has the order  of the $0$ characteristic value less than $k$ and $I_y$ has the order of the $0$ characteristic value less than $s+1.$ But in this case $R_y$ has the order of the $0$ characteristic value less than $k+s,$ i.e. $dim L_0(y) \leq k+s-1.$ Therefore we come to a contradiction with the regularity of the element $a,$ if we suppose that $\overline a$ is not regular. The proof is complete.
\end{proof}

It should be noted that the preimage under the natural homomorphism of a regular element (Cartan subalgebra) is not necessarily regular element (respectively, Cartan subalgebra).

\textbf{Acknowledgments.} \emph{The second and third named authors
would like to acknowledge the hospitality of the "Institut f\"{u}r
Angewandte Mathematik", Universit\"{a}t Bonn (Germany). This work
is supported in part by the DFG 436 USB 113/10/0-1 project
(Germany) and the Fundamental Research Foundation of the
Uzbekistan Academy of Sciences.}

\end{document}